\documentclass[11pt]{article}

\usepackage{amssymb}
\usepackage{graphicx}

\begin{document}

\title{{\bf On causality and closed geodesics of compact Lorentzian manifolds and static spacetimes}}

\author{Miguel S\'anchez
\thanks{ 
Referees' careful reading and comments are acknowledged. Partially supported by MCyT-FEDER Grant n$^o$ MTM 2004-04934-C04-01.
}
\\ Departamento de Geometr\'{\i}a y Topolog\'{\i}a \\ 
Facultad de Ciencias, Universidad de Granada \\
E-18071 Granada, Spain.\\
E-mail: {\ttfamily sanchezm@ugr.es}}

\date{  }
\textwidth 5.5 in
\textheight 8.0 in
\topmargin 0 in
\oddsidemargin 0.7 in
\evensidemargin 0.7 in

\makeatletter

\def\theequation{\thesection.\arabic{equation}}
\@addtoreset{equation}{section}

\newtheorem{defn}{Definition}
\def\thedefn{\thesection.\arabic{defn}}
\@addtoreset{defn}{section}
             
\newtheorem{teor}[defn]{Theorem}
\newtheorem{ejem}[defn]{Example}
\newtheorem{lema}[defn]{Lemma}
\newtheorem{rema}[defn]{Remark}
\newtheorem{coro}[defn]{Corollary}
\newtheorem{prop}[defn]{Proposition}

\makeatother
\font\ddpp=msbm10  at 11 truept 
\def\R{\hbox{\ddpp R}}     
\def\C{\hbox{\ddpp C}}     
\def\L{\hbox{\ddpp L}}    
\def\S{\hbox{\ddpp S}}
\def\Z{\hbox{\ddpp Z}}
\def\Q{\hbox{\ddpp Q}}     

\newcommand{\dps}{d_{\psi}}
\newcommand{\df}{d_{\phi}}

\newcommand{\M}{{\cal M}}
\newcommand{\Mo}{{\cal M}_0}
\newcommand{\be}{\begin{equation}}
\newcommand{\ee}{\end{equation}}
\newcommand{\la}{\Lambda}
\newcommand{\dem}{{\bf Proof: }}
\newcommand{\inte}{\int_{0}^{1}}
\newcommand{\gam}{\gamma}
\newcommand{\eps}{\epsilon}
\newcommand{\<}{\langle}
\renewcommand{\>}{\rangle}
\newcommand{\Om}{\Omega^1}
\renewcommand{\(}{\left(}
\renewcommand{\)}{\right)}
\renewcommand{\[}{\left[}
\renewcommand{\]}{\right]}
\newcommand{\om}{\omega}
\newcommand{\me}{\frac{1}{2}}
\newcommand{\Mt}{\widetilde{\M}}
\newcommand{\cat}{{\mathop{\rm cat}\nolimits}}

\hyphenation{Lo-rent-zian}

\maketitle
\begin{abstract}
\noindent 
Some results related to the causality of compact Lorentzian manifolds are proven: (1) any compact  Lorentzian manifold which admits a timelike conformal vector field is totally vicious, and (2)  a compact Lorentzian manifold covered regularly by a globally hyperbolic spacetime admits a timelike closed geodesic, if some natural topological assumptions (fulfilled, for example, if one of the  conjugacy classes of deck transformations containing a closed timelike curve is finite) hold. As a consequence, any compact Lorentzian manifold conformal to a static spacetime is geodesically connected by causal geodesics, and admits a timelike closed geodesic. \\

\end{abstract} 

\smallskip

\noindent {\em Running title}: Geodesics on compact Lorentzian manifolds.\\

\smallskip





\newpage
 
\section{Introduction}\label{s1}

The aim of this paper is to give some general results on causality and existence of closed timelike geodesics in compact Lorentzian manifolds (Theorems \ref{t1}, \ref{t2}, 
Propositions \ref{p1}, \ref{p2}). Even though these results are interesting in themselves, the principal application will hold for compact static spacetimes (Theorem \ref{t0})\footnote{Any compact spacetime does contain closed timelike curves (CTCs). From a classical relativistic viewpoint, the existence of CTCs is a drawback for a spacetime because of well known paradoxes \cite[p.189]{HE}. Nevertheless, for different reasons there has been a continued interest in spacetimes with CTCs: the existence of CTCs in classical spacetimes such as G\"odel's or the inner part of Kerr's, technical advantages of compactifications, speculations on time-machines and wormholes, quantum interpretations, the recent role of the G\"odel solution as an exact model in string theory (see for example \cite{Ev, Marolf, Vi, Yu})... Nevertheless, our study will remain at a geometrical level.}.

A well-known result by Tipler \cite{Ti} (see also \cite[Theorem 4.15]{BEE}) asserts that any compact Lorentzian manifold, covered regularly by a globally hyperbolic manifold which admits a {\em compact} Cauchy hypersurface, must contain a closed timelike geodesic. This result was extended by Galloway \cite{Ga1}, who also introduced the notion of {\em stable} free t-homotopy class (see also \cite{Ga2}  for related results). Recently, Guediri \cite{GuZ}  has shown that the hypothesis on compactness in Tipler's result cannot be removed, by means of a counterexample (see also \cite{CDB}, \cite{GuZ2}, \cite{GuTr}). Nevertheless, the compactness hypothesis can be replaced by the following assumption  \cite[Theorem 5.1]{GuZ}:  a free t-homotopy class is determined by a central deck transformation $\phi$, i.e., $\phi$ is the unique element in its conjugacy class ${\cal C}$.  Later on, Caponio et al. \cite{CMP} have studied  compact static spacetimes by using some variational results. Essentially, they show that such a spacetime is geodesically connected and, if a free homotopy class is determined by a finite conjugacy class  of deck transformations ${\cal C}$, it contains a closed  geodesic (not necessarily timelike).

\smallskip
\noindent The results in the present article can be summarized as follows.

In Section \ref{s2} some preliminar properties are recalled, with particular emphasis in the geometrical and topological properties of static spacetimes, which will be necessary to apply our results. Especially, Theorem \ref{lll} and Corollary \ref{p0} characterize when the universal covering of a static spacetime is  standard static, and its main properties.

In Section \ref{s3} we prove the following result on the causality of a class of Lorentzian manifolds. Recall that a (time-oriented) Lorentzian manifold $(M,g)$ is called totally vicious if the chronological future and past of any point is the whole manifold, i.e., $I^+(p)=$ $I^-(p) = M, \forall p \in M$.
\begin{teor}\label{t1}
Any compact Lorentzian manifold $(M,g)$ which admits a timelike conformal vector field  is totally vicious.
\end{teor}
The technique of the proof involves some properties of conformal vector fields studied in 
\cite{RS1}. Theorem \ref{t1} will be essential to prove  not only that a compact static manifold is   geodesically connected, but also that any two points can be joined by a timelike geodesic (see Theorem \ref{t0}(1)).

In Section \ref{s4}, we give two extensions of Tipler's result, where the compactness of the 
Cauchy hypersurface is replaced by different assumptions on the group of deck transformations, Propositions \ref{p1}, \ref{p2}. As an immediate consequence of the first one 
(Proposition \ref{p1}), we have the following generalization of Guediri's criterion for the existence of timelike closed geodesics \cite[Section 5]{GuZ}:

\begin{teor} \label{t2}
Let $(M,g)$ be a compact Lorentzian manifold which admits a regular covering $\Pi: \bar M\rightarrow M$ such that $\bar M$ is globally hyperbolic, and let $G$ be the group of deck transformations of $\bar M$. Assume that a 
conjugacy class ${\cal C} \subset G$ satisfies:
\begin{list}{(\alph{enumi})}{\usecounter{enumi} 
\itemsep 0pt \parsep 0pt}
\item It contains a closed timelike curve $\alpha$.
\item ${\cal C}$ is finite.
\end{list}
Then there exists at least one closed timelike geodesic in ${\cal C}$. 
\end{teor}
As a consequence of the second extension (Proposition \ref{p2}), in Section \ref{s5} the results on closed geodesics in \cite{CMP} will be improved by showing that, in a compact static spacetime, a closed {\em timelike} geodesic exists, {\em without any further assumption on the fundamental group}. Moreover, our results will also hold under conformal transformations, because are based 
only on causal and topological properties (in the spirit of \cite{Sa-pams}).
Summing up, we will  prove and discuss:

\begin{teor} \label{t0}
Let $(M,g)$ be a compact static spacetime. Then:

(1) Any pair of points $p,q \in M$ can be joined by means of a timelike geodesic.


(2) Any conjugacy class ${\cal C} \subset G$ which contains a closed timelike curve  contains a closed timelike geodesic too. In particular, there exists at least one closed timelike geodesic in $M$.
\end{teor}

\section{Preliminaries. Static spacetimes}\label{s2}

All Lorentzian manifolds are assumed to be connected, time-oriented (thus, time-orientable), with dimension $n\geq 2$. As usual,  differentiability $C^\infty$ will be assumed, even though, in principle, we only need $C^1$ (geodesics and causality are then well-defined). Our notation and conventions will be standard in Lorentzian Geometry, as in the books \cite{BEE}, \cite{O}, \cite{SW}. A Lorentzian manifold will be called {\em stationary} if it admits a timelike Killing vector field $K$, and {\em static} if, additionally, $K$ is {\em irrotational} (the orthogonal distribution to $K$ is involutive). Standard properties of such manifolds can be seen in 
\cite{SW}, and a survey in \cite{Sa-grec}; for  recent references on the static case 
see\footnote{A brief survey of the static case is carried out in \cite{Sa04}, including an announcement of the results in the present article.} \cite{AU}, \cite{CMP},  \cite{Sa-nonlin}.

The problem of the geodesic connectedness of a Lorentzian manifold has been widely studied recently, specially since Masiello's book \cite{Mas}, which develops a variational viewpoint  (see \cite{Sa-cata} for a survey). Nevertheless, our results on connectedness will rely on a classical theorem by Avez \cite{Av} and Seifert \cite{Se} for causal geodesics: {\em in any globally hyperbolic spacetime, each two causally related points $p, q$ can be joined by a causal geodesic, with length equal to the time-separation (or Lorentzian distance) between $p$ and $q$}.

If the Lorentzian manifold $(M,g)$ is compact, it is well-known that the Euler characteristic of $M$ vanishes. Even though this is not a restriction if the dimension $n$ is odd, it yields a first topological restriction for even $n$; in particular, if $n=2, 4$ then $M$ cannot be simply connected. The condition of stationarity yields new topological obstructions \cite{RS1}; for example, if $(M,g)$ is compact, stationary and $n=3$ then $M$ is a Seifert manifold. Nevertheless, by using Hopf fibration it is not difficult to construct stationary metrics on any odd-dimensional sphere \cite{RS2} (see  \cite{GPR} for further properties). Thus, taking the product of such a stationary 3-sphere by any Riemannian $k$-sphere, with $k>1$, we have:  
{\em there exist simply connected compact stationary manifolds of any dimension $n\geq 5$ and $n=3$}. 

The situation is radically different in the static case. In fact, if $(M,g)$ is a static manifold and $K$ is the corresponding ``static'' (irrotational timelike Killing) vector field, then $K$ is parallel for the conformal metric
\be \label{egstar} g^* = - \frac{1}{g(K,K)} g \ee
(notice that any static spacetime is locally isometric to a standard one $\R\times S$ endowed with a metric as (\ref{egh}) below, with $K$ identifiable to $\partial_t$). Thus, the associated one-form
\be \label{eom} \omega = -g^*(K,\cdot) \ee
is closed and, if $M$ is compact, then it cannot be simply connected. Moreover, the following structural result holds (compare with \cite[Section 3]{CMP}): 

\begin{teor}\label{lll}
Let $(M,g)$ be a 
static manifold with static vector field $K$,  
and $(\bar M,\bar g)$, $\Pi: \bar M \rightarrow M$, $\bar g= \Pi^* g$,
its universal 
Lorentzian covering.  

(1) If the vector field $K$ is  complete 
then $(\bar M,\bar g)$ is  a standard static manifold. More precisely, $\bar M$ is isometric to a product $ \R \times S$ endowed with the metric 
\be \label{egh}
\bar g[(t,x)] = -\beta(x) dt^2 + g_S[x]
\ee
where $g_S$ is a  Riemannian metric on $S$, $dt= \Pi^*\omega$, $K_{\Pi(p)} = d\Pi_p(\partial|_p)$  and, being $\Pi_S: \R \times S \rightarrow S$ the natural projection, $\beta(\Pi_S(p)) = -g(K_{\Pi(p)}, K_{\Pi(p)})$, for all $p \in \bar M$.

(2) If the metric $g$ is (geodesically) complete, 
then  $K$ is complete and the metric $g_S$ in (\ref{egh}) is complete.
\end{teor}
\dem 
(1) Let $\bar K$ be the (complete) vector field on $\bar M$ such that $\Pi_*\bar K = K$, and 
let $\bar \Phi$ be its global flow. As $\bar M$ is simply connected, the closed form 
$\Pi^*\omega$ is exact, i.e., $\Pi^*\omega = dt$ for some function $t: \bar M \rightarrow \R$. 
Fixing $p \in \bar M$  one has $t(\bar \Phi_s(p)) = s + t(p)$ for all $s\in \R$
(use  $dt(K)\equiv 1$). Thus, 
putting $S=t^{-1}(0)$, it is straightforward to check that the required isometry is:
$$ \bar M \rightarrow \R \times S , \quad \quad p \mapsto (t(p), \bar \Phi_{-t(p)}(p)).$$


(2) Let us see  that the vector field $K$ must be complete. Otherwise, for some $p\in M$, a local flow $\Phi$ of $K$ will satisfy that the curve $\lambda \rightarrow \Phi_\lambda(p)$ is well defined for $\lambda\in [0,1)$ but cannot be continuously extended to $\lambda =1$. By using the local decomposition of $M$ as a standard static spacetime, there exists a neighborhood 
$U$ of $p$ isometric to $(-\nu, \nu) \times S_p$, for some $\nu >0$, endowed with a metric as (\ref{egh}). Now, for a small $\mu>0$ ($\mu<\nu \leq 1$) there exists a geodesic 
$$\gamma:[0,1]\rightarrow U \equiv (-\nu, \nu) \times S_p \; \; \mbox{with}  \; \gamma(0)=p (\equiv (0,p)) , 
 \; \gamma(1)=\Phi_{\mu}(p) (\equiv (\mu, p)).$$
For each $\lambda \in [0,1-\mu]$, consider the ({\em complete}) geodesic $\gamma_\lambda$ with initial condition:
$$ \gamma_\lambda(0) = \Phi_\lambda(p) \quad   \quad \quad \mbox{and} \quad \quad  \gamma_\lambda'(0)= d\Phi_\lambda(\gamma'(0)).$$ 
Clearly, for some $\lambda_0>0$ one has:
\be \label{aux}
 \gamma_\lambda(s) = \Phi_\lambda \circ \gamma(s) , \quad \quad \forall s \in [0, 1], \quad
\forall \lambda\in [0,\lambda_0].
\ee
Assume that the domain $V_\lambda$ of $\Phi_\lambda$ ($\Phi_\lambda: V_\lambda \rightarrow M$) is  a maximal neighbourhood of $p$, and define an interval $I\subseteq [0,1-\mu]$ by: $\lambda_0\in [0,1-\mu]$ belongs to $I$ if and only if (\ref{aux}) holds. The result will hold if $I=[0,1-\mu]$ because, in this case, 
$$ \gamma_{1-\mu}(1)=
\Phi_{1-\mu}(\gamma(1))= \Phi_1(p),$$
in contradiction with the inextendibility of $\Phi_\lambda(p)$. As $\gamma([0,1])$ is compact, the interval $I$ is an open subset of $[0,1-\mu]$. Therefore, if $I\neq [0,1-\mu]$, then $I=[0,\lambda_{max})$ for some $0<\lambda_{max}\leq 1$, and the following contradiction would appear. 
Taking into account that the limit of $ \gamma'_\lambda(0)$ when $\lambda \nearrow \lambda_{max}$ is $\gamma'_{\lambda_{max}}(0)$, one has
$$ \gamma_{\lambda_{max}}(s) = \lim_{\lambda\nearrow \lambda_{max}}\Phi_\lambda \circ \gamma(s),
\quad \quad \forall s \in [0,1]
$$
and each integral curve $\lambda \rightarrow \Phi_\lambda(\gamma(s))$ can be continuously extended beyond $\lambda_{max}$, that is, the maximal domain $V_{\lambda_{max}}$ of $\Phi_{\lambda_{max}}$ contains $\gamma(s)$ for all $s\in [0,1]$, and $\lambda_{max}\in I$. 

For the last assertion, recall that any maximal integral manifold $S$ of the kernel of $\omega$ will be complete, because $S$ is totally geodesic in $M$. $\square$

\begin{rema} \label{remarema} {\rm
(1) If $M$ were compact then not only $K$ would be complete but the  static  metric $g$ would be complete too (see \cite{RS1}); thus, Theorem \ref{lll}(2) would be applicable.

(2) A static standard manifold $(\bar M, \bar g)$ as in (\ref{egh}) is globally hyperbolic if 
$g_S$ is complete and $\beta$ behaves at most quadratically at infinity.
[Recall the definition: let $d_R$ be the distance on $S$ canonically associated to the Riemannian metric $g_S$, and assume that, 
 for some fixed $x_0 \in S$ and $ k, k'\in \R$, $p>0$, 
\be \label{quad} \beta(x) \leq k d_R^p(x,x_0) + k'; 
\ee
if (\ref{quad}) holds for $p=2$ (resp. some $p<2$) then $\beta$ is said to behave {\em at most quadratically} (resp {\em subquadratically}) at infinity.] 
In fact, the conformal metric $g_S^*= g_S/\beta$ would be complete too and, thus,  $\bar g^*= \Pi^*(g^*)$ would be globally hyperbolic,  each slice $\{t_0\}\times S$ being
a spacelike Cauchy hypersurface\footnote{Even though a {\em (smooth) spacelike} Cauchy hypersurface exists in any globally hyperbolic spacetime, this is not as trivial as it sounds \cite{BS}.} \cite[Theorem 3.67]{BEE}. As $\bar g$ is globally conformal to $\bar g^*$,    $(\bar M,\bar g)$  is globally hyperbolic too. In particular, this happens for $(\bar M, \bar g)$ in Theorem \ref{lll} if $M$ is compact, 
because $g$ would be complete and Sup$(\beta) <\infty$.

(3) As $dt$ is the pull-back of $\omega$,  any deck transformation $\phi$ of $\bar M$ must preserve $dt$ (i.e., $dt = \phi^*dt = d(t\circ \phi)$), and $t\circ \phi = t + T_\phi$, for some $T_\phi \in \R$, i.e.:
$$
\phi(t,x) = (t+ T_\phi, \Pi_S(\phi(t,x))), \quad \quad \forall (t,x) \in \R\times S,
$$
By deriving partially in both sides with respect to $t$, and taking into account that 
$\phi_*(\partial_t) = \partial_t$, we can write: $\Pi_S(\phi(t,x)) = \phi^S(x)$ (independent of $t$) for some diffeomorphism $\phi^S$ of $S$. 
}\end{rema}
Summing up for the compact case:

\begin{coro} \label{p0}
Let $(M,g)$ be a compact static manifold, and $(\bar M,\bar g)$ its universal Lorentzian covering. Then:

(1) $(\bar M,\bar g)$ is isometric to a globally hyperbolic standard static spacetime $\R \times S$ as in (\ref{egh}), being each slice $\{t_0\}\times S$
a Cauchy hypersurface.

(2) Any deck transformation $\phi: \bar M \rightarrow \bar M$ can be written as
$$
\phi(t,x)= (t+T_\phi , \phi^S(x)),
$$
for some diffeomorphism  $\phi^S$ of $S$ and $T_\phi \in \R$.

\end{coro} 

\begin{ejem} \label{rrr1} {\rm
Notice that $S$ is not necessarily compact. It is especially easy to construct examples in a  torus (recall that any stationary surface is static, because the orthogonal distribution to the timelike Killing vector field $K$ is 1-dimensional; more general examples can be constructed obviously by taking this surface as the fiber of a warped product -or as the base, provided that the warping function is invariant by the flow of $K$). In fact, it is trivial that any flat Lorentzian torus admits a Killing (indeed, parallel) timelike vector field $K$, such that the integral curves of $K^\perp$ (which are isometric to $S$) are not closed. Of course, in this example there are other $K$'s where the curves are closed. But one can also construct a stationary torus with only one independent Killing vector field $K$ such that the integral curves of $K^\perp$ are not closed, as follows. 
Consider $\R^2$, endowed with the Lorentzian metric
$$
g= F(x)(dx \otimes dy + dy \otimes dx) - G(x) dy^2, \quad \quad (F(x) \neq 0, \forall x\in \R ),
$$
where $F, G$ are periodic functions of period 1, and let
$T^2$ be the  Lorentzian torus obtained as the quotient  $\R^2/\Z^2$ (these metrics, as well as those in Remark \ref{r1} below, are particular cases of Lorentzian tori admitting a Killing vector field, studied systematically in \cite{Sa-trans}). The Killing vector field $\partial_y$ projects onto a Killing vector field $K$ on $T^2$. If $G>0$, $K$ is timelike, and $T^2$ is static. By  \cite[Theorem 4.2]{Sa-trans}, if $G' \not\equiv 0$, the metric is not flat and   any other Killing vector field on $T^2$  is a multiple of $K$. Now, recall that the vector field 
$G(x)\partial_x + F(x)\partial_y$ on $\R^2$, projects onto a non-vanishing 
vector field $W$ on $T^2$ orthogonal to $K$. Finally, it is easy to check that if
$$
\int_0^1 \frac{F}{G}(x) dx 
$$
is not rational, then the integral curves of $W$ are not closed, as required.
}\end{ejem}

\begin{rema} \label{rrr2} {\rm 
Some additional information on 
$(M,g)$ in Corollary \ref{p0} can be obtained, in comparison with the general results in 
\cite{Ze}. Recall that the Levi-Civita connection $\nabla^*$ of $g^*$ in (\ref{egstar}) 
is Riemannian, i.e., the Riemanian metric 
$g^*_R(A,B)= g^*(A,B) -2g^*(A,K)g^*(B,K)$ has the same Levi-Civita connection that $\nabla^*$.
Then, deck transformations for $\bar M$ are also isometries for both, $\bar g^*= \Pi^*(g^*)$ and $\bar g^*_R= \Pi^*(g^*_R)$. Now, write $(\bar M, \bar g^*)$ as a semi-Riemannian product $\L^k\times N$ where $\L^k, k\geq 1$, is a $k$-dimensional Lorentz Minkowski spacetime ($\partial_t$ will be chosen to project on $K$), and $N$ is a Riemannian manifold with no further decomposition as a Riemannian product ($N\neq N'\times \R$). Thus, any deck transformation $\phi$ of $\bar M$ can be written as a composition $\phi_1\circ \phi_2$, where $\phi_1$ is an isometry of $N$
and $\phi_2$ is an isometry of $\L^k$ (and $\R^k$) which preserves $\partial_t$ [i.e., $\phi_2$ can be identified to an element of the  semi-direct product $O(k-1,\R ) \times \R^k= 
(O(k,\R) \cap O^{\uparrow}_1(k,\R))\times \R^k$].
}\end{rema}
 
\section{Connecting timelike curves}\label{s3}

Notice that a timelike conformal vector $K$ for $g$ is Killing and unitary for the conformal metric $g^*=  -(1/g(K,K))g, $
(see \cite[Lemma 2.1]{Sa-trans}). 
Thus, Theorem \ref{t1} is equivalent to:

\begin{prop}
Any compact Lorentzian manifold $(M,g)$ admitting a Killing vector field $K$ with $g(K,K)=-1$ is totally vicious.
\end{prop}
\dem It is not difficult to prove that a Lorentzian manifold is totally vicious if and only if for every point $p$ there exists a closed timelike curve through $p$ (see  \cite[Proposition 2.2]{Ma}). Thus, it is enough to show  that $(M,g)$ admits a timelike vector field $X \in \Gamma(TM)$ with closed integral curves. 
Consider the auxiliary Riemannian metric
$
g_R(A,B)= g(A,B) - 2g(A,K) g(B,K)
$
for all $A,B \in \Gamma(TM)$, which will have a compact isometry group Iso$(M,g_R)$. A straightforward computation shows that $K$ is also a Killing vector field for $g_R$ and, thus, its one-parameter group $G$ has a compact closure $\bar G$ in Iso$(M,g_R)$. As $G$ is abelian, $\bar G$ is  abelian too and, thus, isomorphic to a $k$-torus, for some $k\geq 1$. Therefore, there is a sequence of one-parameter subgroups $\{ G_m\}$ diffeomorphic to circles, whose associated sequence of $g_R$--Killing vector fields $\{ X_m\}$ converges to $K$ 
(i.e., $\lim_{m\rightarrow \infty} {\rm Max}_{p\in M} \; g_R(X_m(p)-K(p), X_m(p)-K(p)) =0$). 
Thus, for some $m_0$ sufficiently large, 
$X_{m_0}$ is timelike, and we can choose $X=X_{m_0}$. $\square$

\begin{rema} \label{r1} {\em
Total viciousness may {\em not} hold if the compact Lorentzian manifold $(M,g)$ is assumed to admit a Killing vector field $K$ which is only {\em causal}. In fact, it is not difficult to construct counterexamples among Lorentzian tori $\R^2/\Z^2$ obtained as a quotient of $\R^2$ endowed with the metric $g = \sin (\psi(x)) (dx^2-dy^2)  + 2 cos (\psi(x)) dxdy $ 
for suitable functions $\psi(x)$ of period 1 (see Fig. 1). 
Recall that, in these counterexamples, $K$ is also irrotational.

}\end{rema}

\vspace{2cm}
\begin{center}

{\bf Fig. 1 to be inserted here.}
\end{center}

\vspace{2cm}

\section{Closed timelike geodesics}\label{s4}

Let $(M,g)$ be a compact Lorentzian manifold and $\Pi: \bar M\rightarrow M$ a regular covering endowed with the pullback metric $\bar g= \Pi^*g$. Assume that $(\bar M,\bar g)$ is globally hyperbolic and, thus, topologically $\bar M = \R \times S$, where $S$ is a Cauchy hypersurface. Let $d: \bar M\times \bar M \rightarrow [0,\infty)$ be the Lorentzian time-separation (or Lorentzian distance)  on $\bar M$. Recall that, because of global hyperbolicity, 
$d$ is  continuous and finite--valued. 
For each deck transformation $\phi \in G$, consider the function
\be \label{edfi} \df: \bar M \rightarrow [0, \infty), \quad \quad \df(p) = d(p, \phi(p)), \quad \forall p \in \bar M.\ee
\begin{lema} \label{l1}
If $\df$  restricted to $S$ admits a relative maximum $p_0$ with $\df(p_0)>0$ then there exists a timelike closed geodesic in $(M,g)$.
\end{lema}
\dem  The result is a consequence of Tipler's technique. In fact, as $\bar M$ is globally hyperbolic and $\df(p_0)>0$ (i.e., $\phi(p_0) \in I^+(p_0)$) the Avez-Seifert result yields a   timelike geodesic $\bar \gamma:[0,1]\rightarrow \bar M$ from $p_0$ to $\phi(p_0)$ which is maximizing, i.e., length$(\bar \gamma )= \df(p_0)$. Then, $\gamma= \Pi\circ \bar \gamma$ is the required geodesic (otherwise, $\gamma'(0) \neq \gamma'(1)$ and $\gamma$ could be modified in any arbitrarily smooth neighborhood of $\gamma(0) (=\gamma(1))$ to obtain a strictly longer closed timelike curve, which contradicts the condition of relative maximum of $p_0$). $\square$.

\smallskip

\noindent Let ${\cal C}$ be a conjugacy class of the group $G$ of deck transformations of $\bar M$. Even though, in general, a closed curve $\alpha$ does not determine any deck transformation,    $\alpha$ does determine a conjugacy class of $G$ and, thus, to assert that ${\cal C}$ contains $\alpha$  makes sense. Even more, in the special case of closed timelike curves, if $\gamma_1$ and $\gamma_2$ are two freely t-homotopic closed curves (in the sense of \cite{Ga1}, i.e., freely homotopic through timelike curves) with base points $x_1$, $x_2$ on each one, and if $\bar x_1$, $\bar x_2$ are two points on $\bar M$ over $x_1$, $x_2$, resp., both belonging to the same Cauchy hypersurface, then $\gamma_1$ and $\gamma_2$ determine the same deck transformation. In particular, given a free t-homotopy class $\tilde {\cal C}$, each Cauchy hypersurface determines a unique deck transformation in the conjugacy class representing $\tilde {\cal C}$.

\begin{prop} \label{p1}
Let $(M,g)$ be a compact Lorentzian manifold which admits  a globally hyperbolic manifold $(\bar M, \bar g)$ as a regular covering. Assume that a conjugacy class of deck transformations ${\cal C}\subset G$ satisfies:

(a) ${\cal C}$ contains a closed timelike curve $\alpha$.

(b) For some (and then for any) compact subset $K \subset \bar M$ such that $M \subset \Pi (K)$ all the restricted functions  $\df|_K$, $\phi \in {\cal C}$ are null, except for a finite number $\phi_1, \dots \phi_j \in {\cal C}$. 

\smallskip

Then there exists a closed timelike geodesic in ${\cal C}$.  

\end{prop}
\dem Fix a Cauchy hypersurface, and let $\phi_0$ be the unique 
deck transformation in ${\cal C}$ determined by the t-homotopy class of $\alpha$. From Lemma \ref{l1}, it is enough that $d_{\phi_0}$ attains an absolute maximum on $\bar M$. Recall first that, for any deck transformation 
$\psi \in G$:
\be \label{e1}
 d_{\phi_0}(\psi(p)) = d(\psi(p), \phi_0 (\psi((p))) = d(p, \psi^{-1}\circ \phi_0 \circ \psi(p)) =
d_{\psi^{-1}\phi_0\psi}(p).
\ee
On the other hand, taking into account that $\Pi(K)= M$ and the covering is regular:
\be \label{e2}
\{ d_{\phi_0}(p): p\in \bar M\} = \{ d_{\phi_0}(\psi(p)): p\in K, \psi \in G\}.
\ee
As only finitely many conjugate $\psi^{-1}\circ \phi_0\circ \psi \in {\cal C}$ are non identically zero on $K$, the supremum in (\ref{e1}) is then equal to the maximum of 
$$
\{d_{\phi_i}(p_i), \quad \quad i= 1, \dots j\},
$$
attained at some index $i_0$, where each $p_i$ is the maximum of the non-null function $d_{\phi_i}|_K$. 
Therefore, 
by  (\ref{e1}) the absolute maximum of $d_{\phi_0}$ is attained at  $\psi_{i_0}(p_{i_0})$, where $\phi_0 = 
\psi_{i_0} \circ \phi_{i_0}\circ \psi_{i_0}^{-1}$. $\square$

\begin{rema} \label{rctc}
{\em
Assumption (a) is natural and not too restrictive, because any compact Lorentzian manifold admits a closed timelike curve (see for example \cite[Lemma 14.10]{O}). Then, the curve $\alpha$  determines a free t-homotopy class where, in fact, a closed timelike geodesic should appear. Of course, assumption (b) is always satisfied if ${\cal C}$ is finite, yielding 
Theorem \ref{t2}  (compare with Sections 4 and 5 in \cite{GuZ}). Nevertheless, the more general assertion in (b) will be needed to prove Proposition \ref{p2}, as we will see below (see Step 2 of the proof of this proposition). 
}\end{rema}


\noindent In Lemma \ref{l1} we saw that the existence of a (relative) maximum of $\df$ on $S$ was enough to obtain a closed timelike geodesic, but in the proof of Proposition \ref{p1} we ensured  a stronger property, the existence of a maximum on all $\bar M$. The reason relies in the lack of good technical properties of $\phi|_S$ 
when $\phi$ is changed by other element of ${\cal C}$.   
Nevertheless, there are interesting spacetimes, as the static ones, where these technical properties (namely, (\ref{erec}) below) occur. In this case, assumption (b) in Proposition \ref{p1} can be dropped:
\begin{prop} \label{p2}
Let $(M,g)$ be a compact Lorentzian manifold which admits a  globally hyperbolic 
manifold $(\bar M, \bar g)$ as a regular covering.
Assume that any deck transformation $\phi \in G$ of $\bar M = \R \times S$ can be written as 
\be \label{erec}
\phi(t,x)= (t+T_\phi , \phi^S(x)),
\ee
for some homeomorphism  $\phi^S$ of $S$ and $T_\phi \in \R$.
 
If the conjugacy class ${\cal C}$  contains a closed timelike curve $\alpha$, then
there exists a closed timelike geodesic in ${\cal C}$.  

\end{prop}
\dem
It is straightforward to check that the value of $T_\phi$ is equal for all the elements in a same conjugacy class; thus, we can write:
\be \label{erec2}
\phi(t,x)= (t+T, \phi^S(x)), \quad \quad \forall \phi \in {\cal C}.
\ee
Let $\phi_0 \in {\cal C}$ be the deck transformation determined by $\alpha$ and $S\equiv \{0\} \times S$, and let $d_S$ be the restriction of  
function $d_{\phi_0}$ to $S$, i.e.: 
\be \label{edese} d_S: S \rightarrow [0, \infty), \quad \quad d_S(x) = d((0,x),(T,\phi_0^S(x)), \quad \forall x \in S. \ee
From Lemma \ref{l1}, it is enough to show that $d_S$ attains a maximum. The proof consists of the following three steps: 

{\em Step 1}: To choose a compact subset $K^S\subset S $ such that:
$$ S = \cup_{\psi \in G} \; \psi^S(K^S). $$
$K^S$ can be chosen as follows. Consider a compact subset $K\subset \R\times S$ such that $\Pi(K)=M$.  
Clearly, $K$ can be chosen as a product $K=[0,\bar T]\times K^S$, and 
$$ \bar M = \cup_{\psi \in G} \; \psi(K) = \cup_{\psi \in G} \; [T_{\psi},T_{\psi}+\bar T]\times \psi^S(K^S).$$
Thus, if $\Pi_S: \R \times S \rightarrow S$ is the natural projection:
$$S= \Pi_S(\bar M) = \cup_{\psi \in G} \; \psi^S(K^S).$$

{\em Step 2}: To show that, fixed $T$ in (\ref{erec2}), there are only finitely many deck transformations $\phi\in {\cal C}$ satisfying 
\be \label{j+}
 (T,\phi^S(K^S)) \cap J^+(0,K^S) \neq \emptyset 
\ee
(i.e., all the $d_{\phi}$'s in (\ref{edfi}) restricted to $(0,K^S)$ are null for $\phi \in {\cal C}$, except for a finite subset of $\phi$'s).  
Notice first that $(T,S) \cap J^+(0,K^S)$ is compact \cite[Lemma 3.1]{Sa-pams}, as well as the following general result: given two compact subsets $K_1, K_2 \subset \bar M$ the set
$$
\{\phi \in {\cal C}: K_1 \cap \phi(K_2) \neq \emptyset \}
$$
is finite (in fact, the result holds even if $\phi$ is allowed to vary in all $G$, see 
\cite[Lemma 3.2]{Sa-pams}). Then, one has just to apply this result to 
$K_1 = (T,S) \cap J^+(0,K^S)$ and $K_2=(0,K^S)$.

{\em Step 3}: Reasoning as in the proof of Proposition \ref{p1}, to prove that $d_S$ in (\ref{edese}) attains a maximum. In fact, if $\phi_1, \dots \phi_j \in {\cal C}$ are the only deck transformations such that $\phi_1^S, \dots \phi_j^S$ satisfy (\ref{j+}), one has:
$$
{\rm Sup}\{d_S(p): p \in S\} = {\rm Sup}\{d_{\psi^{-1}\phi_0\psi}(p): p \in K^S, \psi \in G\}
$$
$$= {\rm Max}\{d_{\phi_i}(p): p\in K^S, i=1,\dots j\}
$$
(in the first equality Step 1 and (\ref{e1}) are used, and in the second, Step 2).
$\square$

\section{Application to  static spacetimes} \label{s5}

{\it Proof of Theorem \ref{t0}.} For (1), recall that Theorem \ref{t1} ensures the existence of a timelike curve $\alpha$ from $p$ to $q$. Lifting $\alpha$ to a curve $\bar \alpha$ in the universal covering $(\bar M,\bar g)$ and using that this is globally hyperbolic (Corollary \ref{p0}), Avez-Seifert result yields a timelike geodesic $\bar \gamma$ connecting the endpoints of $\bar \alpha$. Thus, the required geodesic is $\gamma = \Pi\circ \bar \gamma$.

The part (2) is obvious from 
Proposition \ref{p2} and Corollary \ref{p0}, plus Remark \ref{rctc}. $\square$

\begin{rema} \label{rf} {\rm 
(1) We saw in Remark \ref{rrr2} that the affine connection $\nabla^*$ on $M$ associated to the conformal metric $g^*$ in (\ref{egstar}) is Riemannian and, thus, $g^*$ is geodesically connected and admits closed geodesics. Nevertheless, this does not imply directly that $g$
also satisfies these two properties because, as far as we know, such properties are not conformally invariant (even on compact manifolds). But, as a clear difference with the techniques in \cite{CMP}, all the properties we have used to prove Theorem \ref{t0} are explicitly conformally invariant (say, hypotheses as {\em (a), (b)} in Proposition \ref{p1} holds for $g$ if and only if hold for any conformal $g^*$, because both metrics have equal timelike vectors and relations of causality), and thus:
\begin{quote}
If $(M,g)$ is a compact static spacetime and $\Omega: M \rightarrow (0,\infty) $ is any 
function, then the conformal metric $g^*= \Omega g$ also satisfies both conclusions (1) and (2) in Theorem \ref{t0}.
\end{quote}
In particular, this is applicable to warped products $(B\times_f F, g=g_B+f^2g_F)$ where one of the factors, say, the base $(B,g_B)$ is a compact static manifold, and the other $(F,g_F)$ a compact Riemannian manifold. In order to check if a conjugacy class type 
${\cal C}_B\times {\cal C}_F$ contains a closed timelike geodesic:
(i) check that ${\cal C}_B$ contains  closed timelike curves, (ii) compute the maximum $L$ of the lengths of CTCs for the conformal metric $g_B^*= g_B/f^2$ (this is equal to the $g_B^*$-length of a maximizing closed timelike $g_B^*$-geodesic in ${\cal C}_B$), (iii) compute the minimum $l$ of the $g_F$-lengths for curves in ${\cal C}_F$ (equal to the length of a minimizing closed $g_F$-geodesic in $F$), and (iv) ${\cal C}_B\times {\cal C}_F$ admits a closed timelike geodesic if and only if $l<L$. 

(2) Essentially, Theorem \ref{t0} improves widely the corresponding results in \cite{CMP} (for example, Corollaries 4.6, 4.7 and 4.8 in \cite{CMP} are particular cases). As suggested by the authors of this reference, an interesting open question would be to determine which conclusions of Theorem \ref{t0} hold if $(M,g)$ is just stationary. Notice that, even though many interesting stationary compact manifolds will satisfy the assumptions of Proposition \ref{p2} (see for example \cite{Sa-pams}), there are others which {\em do not} satisfy them, as the simply connected ones in Section \ref{s2}. On the other hand, recall that the question whether a compact Lorentzian manifold admits a closed (non-necessarily causal) geodesic 
\cite{Ga2} remains open, as far as we know.

(3) In the non--compact case, the authors of \cite{CMP} use the following result (essentially contained in \cite{GP}, see also \cite{Mas}): {\it a standard static spacetime $\R\times S$ as in (\ref{egh}) with $g_S$ complete and $\beta$ subquadratic is geodesically connected}. Recall that, in this case,  the spacetime is globally hyperbolic too (Remark \ref{remarema}(2)).
From the results in 
\cite{FS}, chosen $\epsilon>0$, there exist counterexamples to geodesic connectedness even if inequality (\ref{quad}) holds with $p=2+\epsilon$. Thus, the quadratic case $p=2$ becomes critical for geodesic connectedness. Nevertheless,  even  in this case it is  possible to prove geodesic connectedness \cite{BCFS}.
}\end{rema}

\newpage

\vspace{6cm}
CAPTION FIG. 1:

\begin{center}
Fig. 1. The Killing vector field $K= \partial_y$ is causal, but the quotient torus is not totally vicious: $I^+(1/3,0) (= I^+(2/3,0)) = ]1/3,2/3[ \times S^1$.
\end{center}

\newpage

\includegraphics[width=15cm]{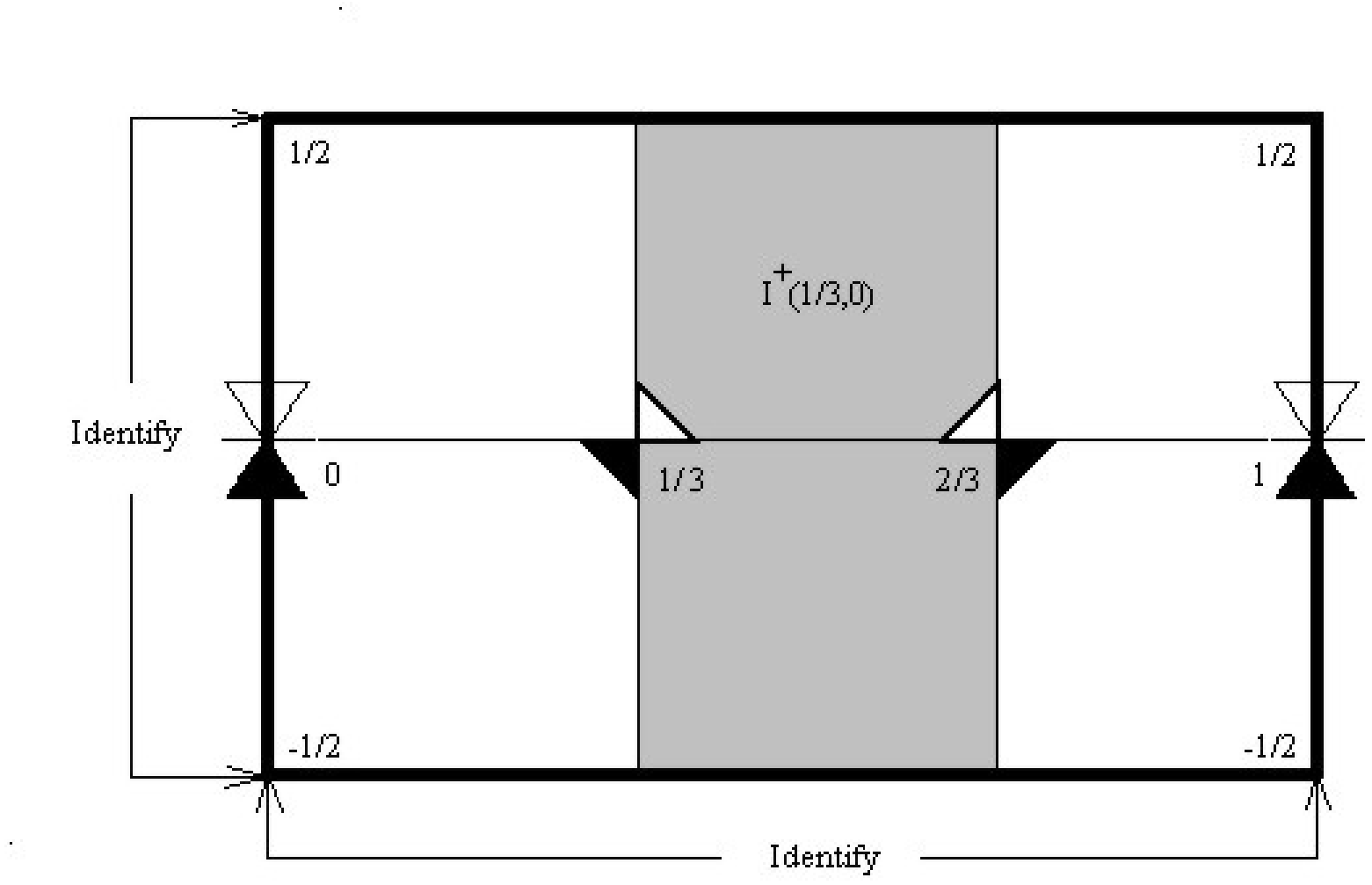}

\end{document}